# Geometric Structures and Differential Operators on Manifolds Having Super tangent Bundle

Naser Boroojerdian

Department of Mathematics and Computer Science,
Amirkabir University of Technology,

Tehran, Iran

**Abstract**

In this paper, we introduce the notion of a super tangent bundle of a manifold, and extend the basic notions of differential geometry such as differential forms, exterior derivation, connection, metric and divergence on manifolds that equipped with a super tangent bundle.

# Content

# 1. Tensors on super spaces

The aim of this paper is to provide some mathematical apparatus needed for calculus in super structures on manifolds. Most of the concepts of this section are found in [2][3]. In this paper, all vector spaces are real or complex finite dimensional. A vector space $V$ with a given decomposition $V = V_0 \oplus V_1$ is called a super vector space. Elements of $V_0$ and $V_1$ are called homogeneous elements of $V$. The parity of a homogeneous element $x$ is denoted by $|x|$, and it is defined as follows.

$$|x| = \begin{cases} 0 & x \in V_0 \\ 1 & x \in V_1 \end{cases}$$

For simplicity, we denote $(-1)^{|x|}$ by $(-1)^x$. In other words, when a homogeneous element of a super space, like $x$, appears on the exponent of $(-1)$, it means $|x|$. When we use the parity of a vector, it is assumed that the vector is homogeneous.

Super structure of $V$ endows some super structures on the spaces of tensors on $V$. For example, $V^* = V_0^* \oplus V_1^*$, so $V^*$ is also a super space. For simplicity, we speak only about covariant tensors, but the same things also hold for contravariant tensors. Elements of $\otimes^k V^*$ are covariant tensors on $V$ of order $k$. Covariant tensors on $V$, can be considered as multilinear maps on $V$. For example, a $k$-linear map $\widetilde{\omega}: V \times \cdots \times V \to \mathbb{R}$ (for real spaces) is a covariant tensors of order $k$. We can consider $\mathbb{R}$ as a super space that its odd part is $\{0\}$. If $\widetilde{\omega}$ preserves parties of vectors, (i.e. $|\widetilde{\omega}(v_1, \cdots, v_k)| = |v_1| + \cdots + |v_k|$) it is called even, and if $\omega$ reverse parties of vectors, (i.e. $|\widetilde{\omega}(v_1, \cdots, v_k)| = |v_1| + \cdots + |v_k| + 1$) it is called odd. Clearly, for $\widetilde{\omega} \in \otimes^k V^*$, it is even if and only if for homogenous vectors $v_1, \cdots, v_k$, if the number of odd vectors is odd, then $\widetilde{\omega}(v_1, \cdots, v_k) = 0$, and $\widetilde{\omega}$ is odd if and only if, if the number of odd vectors is even, then $\widetilde{\omega}(v_1, \cdots, v_k) = 0$.

Here, we interest in the alternating tensors. In the context of super spaces, a covariant tensor $\widetilde{\omega}$ of order $k$, is called (super)alternating, if for homogenous vectors $v_1, \cdots, v_k$, and indices $i = 1, \cdots, k-1$

$$\widetilde{\omega}(v_1, \cdots, v_i, v_{i+1}, \cdots, v_k) = -(-1)^{v_i v_{i+1}} \widetilde{\omega}(v_1, \cdots, v_{i+1}, v_i, \cdots, v_k) \quad (1)$$

The formula given above, for two arbitrary indices $i$ and $j$ is equivalent to the following.

$$\widetilde{\omega}(v_1, \cdots, v_i, \cdots, v_j, \cdots, v_k) = -(-1)^{v_i v_j + (v_i + v_j)(v_{i+1} + \cdots + v_{j-1})} \widetilde{\omega}(v_1, \cdots, v_j, \cdots, v_i, \cdots, v_k) \quad (2)$$

The set of all alternating tensors of order $k$ on $V$ is a vector space which is denoted by $\hat{A}^k(V)$. Alternating tensors are symmetric with respect to odd vectors and are alternating with respect to even vectors. In fact:

$$\hat{A}^k(V) = \sum_{l+s=k} (\wedge^s V_0^*) \otimes (V^l V_1^*) \quad (3)$$

The characteristic property of $\widetilde{\omega} \in (\wedge^s V_0^*) \otimes (V^l V_1^*)$ is that for homogenous vectors $v_1, \cdots, v_{s+l}$, if the number of even vectors is not $s$ (or equivalently the number of odd vectors is not $l$), then $\widetilde{\omega}(v_1, \cdots, v_{s+l}) = 0$. In this case, we call $\widetilde{\omega}$ has the pair order $(s, l)$ so, the total order of $\widetilde{\omega}$ is $s + l$. For $\omega \in \wedge^s V_0^*$ and $\xi \in V^l V_1^*$ we construct $\omega \otimes \xi$ as an element of $(\wedge^s V_0^*) \otimes (V^l V_1^*)$ as follows. For homogenous vectors $v_1, \cdots, v_{s+l}$ that $v_1, \cdots, v_s$ are even and $v_{s+1}, \cdots, v_{s+l}$ are odd, define:

$$(\omega \otimes \xi)(v_1, \cdots, v_{s+l}) = \omega(v_1, \cdots, v_s) \xi(v_{s+1}, \cdots, v_{s+l}) \quad (4)$$

We can extend uniquely this operator to all vectors such that it becomes an element of $(\wedge^s V_0^*) \otimes (V^l V_1^*)$. These tensors are called simple tensors. It can be shown that all alternating tensors are sum of some simple tensors.

For alternating tensors of the pair order $(s, l)$, we can define two parities. The first parity is $s \pmod{2}$, and the second parity is $l \pmod{2}$. The first parity of $\widetilde{\omega}$ is denoted by $|\widetilde{\omega}|$ and the second parity of $\widetilde{\omega}$ is denoted by $|\widetilde{\omega}|'$. Also we denote $(-1)^{|\widetilde{\omega}|}$ by $(-1)^{\widetilde{\omega}}$, and denote $(-1)^{|\widetilde{\omega}|'}$ by $(-1)^{\widetilde{\omega}'}$.

The second parity of a tensor shows that as an operator, it is even or odd. It is sufficient to check this assertion for simple tensors.

Exterior product of (super) alternating tensors can naturally be defined for simple tensors and extend to all alternating tensors. For simple alternating tensors $\omega \otimes \xi$ and $\varphi \otimes \eta$ set:

$$(\omega \otimes \xi) \wedge (\varphi \otimes \eta) = (\omega \wedge \varphi) \otimes (\xi \vee \eta) \quad (5)$$

If the pair order of $\widetilde{\omega}$ is $(s, l)$, and the pair order of $\widetilde{\varphi}$ is $(s', l')$, then the pair order of $\widetilde{\omega} \wedge \widetilde{\varphi}$ is $(s + s', l + l')$. So, this product preserves the first and second parity. This product is associative and super commutative with respect to the first parity.

$$\tilde{\omega} \wedge \tilde{\varphi} = (-1)^{\tilde{\omega}\tilde{\varphi}} \tilde{\varphi} \wedge \tilde{\omega} \qquad (6)$$

The set of all alternating tensors of all orders on $V$, is denoted by $\hat{A}(V)$. This is a super algebra with respect to the exterior product and either of parities.

**Vector valued tensors:** The concept of tensors on the super spaces can be extended to vector valued tensors. If $W$ is some other super space (maybe $V$ itself), a $k$-linear map $\boldsymbol{\omega}: V \times \cdots \times V \to W$ is called a $W$-valued covariant tensors of order $k$ on $V$. If $\boldsymbol{\omega}$ preserve parities of vectors, it is called even, and if $\boldsymbol{\omega}$ reverse parities of vectors, $\boldsymbol{\omega}$ it is called odd.

Alternating vector valued tensors are defined similarly and the set of all $W$-valued alternating tensors of order $k$ on $V$ is a vector space that is denoted by $\hat{A}^k(V, W)$. In fact $\hat{A}^k(V, W) = \hat{A}^k(V) \otimes W$. For $\tilde{\omega} \in \hat{A}^k(V)$ and $w \in W$ define $\tilde{\omega} \otimes w$ as an element of $\hat{A}^k(V, W)$ as follows.

$$(\tilde{\omega} \otimes w)(v_1, \cdots, v_k) = \tilde{\omega}(v_1, \cdots, v_k) w \qquad (7)$$

These tensors are called simple vector valued tensors and generate $\hat{A}^k(V, W)$. Ordinary tensors are number valued tensors, because for $W = \mathbb{R}$ (odd part of $W$ is $\{0\}$), $W$-valued tensors are ordinary tensors. A vector valued tensor $\boldsymbol{\omega}$ of order $k$, is said to have the pair order $(s, l)$, whenever $s + l = k$ and for homogenous vectors $v_1, \cdots, v_k$, if the number of even vectors is not $s$ (or equivalently the number of odd vectors is not $l$), then $\boldsymbol{\omega}(v_1, \cdots, v_{s+l}) = 0$. For example, if $\tilde{\omega}$ is an ordinary tensor of the pair order $(s, l)$ then $\tilde{\omega} \otimes w$ also has the pair order $(s, l)$. For a simple vector valued tensor such as $\tilde{\omega} \otimes w$, the first and second parity of it, is defined as follows.

$$|\tilde{\omega} \otimes w| = |\tilde{\omega}| \quad , \quad |\tilde{\omega} \otimes w|' = |\tilde{\omega}|' + |w| \qquad (8)$$

This definition implies that a vector valued tensor $\boldsymbol{\omega}$ as an operator is even(odd) if and only if $|\boldsymbol{\omega}|' = 0 (|\boldsymbol{\omega}|' = 1)$. The set of all $W$-valued alternating tensors on $V$ is denoted by $\hat{A}(V, W)$ and it is equal to $\hat{A}(V) \otimes W$.

If $W_1, W_2, W_3$ are three super space and $(u, v) \mapsto uv$ is a bilinear map $W_1 \times W_2 \to W_3$, that preserve parities (i.e. $|uv| = |u| + |v|$), associated to it, can be defined an exterior product between $W_1$-valued and $W_2$-valued tensors to produce a $W_3$-valued tensor. For simple tensors $\tilde{\omega}\otimes u$ and $\tilde{\varphi}\otimes v$, the exterior product of them is defined as follows.

$$(\tilde{\omega}\otimes u) \wedge (\tilde{\varphi}\otimes v) = (\tilde{\omega} \wedge \tilde{\varphi})\otimes(uv) \tag{9}$$

It can easily be shown that this product preserve both parities. For example, the product of scalars to vectors produces an exterior product between numerical valued tensors and vector valued tensors. For $\tilde{\varphi} \in \hat{A}^k(V)$ and a simple vector valued tensor $\tilde{\omega}\otimes u$ we have:

$$\tilde{\varphi} \wedge (\tilde{\omega}\otimes u) = (\tilde{\varphi} \wedge \tilde{\omega})\otimes u \tag{10}$$

**Insertion operators:** For $\tilde{\omega} \in \hat{A}^k(V)$ and $v \in V$, insertion of $v$ in $\tilde{\omega}$ is denoted by $i_v\tilde{\omega}$ which is an element of $\hat{A}^{k-1}(V)$. The tensor $i_v\tilde{\omega}$ as an operator, is defined by the following formula.

$$(i_v\tilde{\omega})(v_1,\cdots,v_{k-1}) = (-1)^{\tilde{\omega}v}\tilde{\omega}(v,v_1,\cdots,v_{k-1}) \tag{11}$$

For the case $k = 0$, it is defined $i_v\tilde{\omega} = 0$. If $\tilde{\omega}$ has the pair order $(s,l)$ and $v$ is even, then the pair order of $i_v\tilde{\omega}$ is $(s-1,l)$. And if $v$ is odd, then the pair order of $i_v\tilde{\omega}$ is $(s,l-1)$. So, we can deduce that:

$$|i_v\tilde{\omega}|' = |\tilde{\omega}|' + |v| \quad , \quad |i_v\tilde{\omega}| = |\tilde{\omega}| + |v| + 1 \tag{12}$$

For a simple tensor $\omega\otimes\xi$, if $v$ is even, then $i_v(\omega\otimes\xi) = (i_v\omega)\otimes\xi$ and if $v$ is odd, then $i_v(\omega\otimes\xi) = \omega\otimes i_v\xi$. The insertion operator, with respect to the exterior product has the following property.

$$i_v(\tilde{\omega} \wedge \tilde{\varphi}) = (i_v\tilde{\omega}) \wedge \tilde{\varphi} + (-1)^{\tilde{\omega}(v+1)}\tilde{\omega} \wedge i_v\tilde{\varphi} \tag{13}$$

If we consider $\hat{A}(V)$ as a super algebra with respect to the first parity, then the parity of the operator $i_v$ is $|v| + 1$ and this operator is a super derivations on $\hat{A}(V)$.

## 2. Supper inner product on the super spaces

First, we remind some basic facts about symmetric and symplectic products on ordinary vector spaces. Suppose that $W$ is an $n$ dimensional real vector space and $<,>$ is a nondegenerate symmetric or symplectic product on it. This product induce the isomorphism $\#: W \to W^*$ defined by $u^{\#}(v) = <u,v>$. The inverse isomorphism is denoted by $\alpha \mapsto \alpha^b$, in fact, $<\alpha^b, v> = \alpha(v)$. By this isomorphism the product on $W$ is transferred to $W^*$.

The product on $W$ can be extended to the tensor spaces on $W$. For example, if $v_1, \cdots, v_k, u_1, \cdots, u_k \in W$, then the product of $v_1 \wedge \cdots \wedge v_k$ and $u_1 \wedge \cdots \wedge u_k$ is defined as follows.

$$<v_1 \wedge \cdots \wedge v_k, u_1 \wedge \cdots \wedge u_k> = \det(<v_i, u_j>) \tag{14}$$

For different $k$ and $l$, it is defined $\wedge^k W^*$ and $\wedge^l W^*$ to be orthogonal. If the product on $W$ is symmetric, then the product on all $\wedge^k W^*$,s are symmetric, but if the product is symplectic then for $\omega, \varphi \in \wedge^k W^*$ we have $<\omega, \varphi> = (-1)^k <\varphi, \omega>$.

The product on $W$ can also be extended to a product on $\vee^l W$. For vectors $v_1, \cdots, v_l, u_1, \cdots, u_l \in W$, the product of $v_1 \vee \cdots \vee v_l$ and $u_1 \vee \cdots \vee u_l$ is defined as follows.

$$\begin{aligned}<v_1 \vee \cdots \vee v_l, u_1 \vee \cdots \vee u_l> &= \text{perm}(<v_i, u_j>) \\ &= \sum_\sigma <v_1, u_{\sigma(1)}> \cdots <v_l, u_{\sigma(l)}>\end{aligned} \tag{15}$$

For different $k$ and $l$, it is defined $\vee^k W^*$ and $\vee^l W^*$ to be orthogonal. If the product on $W$ is symmetric, then the product on all $\vee^l W^*$,s are symmetric, but if the product is symplectic then for $\xi, \eta \in \vee^l W^*$ we have $<\xi, \eta> = (-1)^l <\eta, \xi>$.

**Theorem (2.1):** Suppose $<,>$ is a product (symmetric or symplectic) on the vector space $W$. For $v \in W$ and $\omega \in \wedge^{k+1} W^*$ and $\varphi \in \wedge^k W^*$ and $\xi \in \vee^{l+1} W^*$ and $\eta \in \vee^l W^*$ the following equalities hold.

$$\begin{aligned}<i_v \omega, \varphi> &= <\omega, v^{\#} \wedge \varphi> \\ <i_v \xi, \eta> &= <\xi, v^{\#} \vee \varphi>\end{aligned} \tag{16}$$

**Proof:** It can easily be shown that for $v_1, \cdots, v_l \in W$ and $\alpha_1, \cdots, \alpha_l \in W^*$ the following relations hold.

$$\begin{aligned}(\alpha_1 \wedge \cdots \wedge \alpha_l)(v_1, \cdots, v_l) &= <\alpha_1^b \wedge \cdots \wedge \alpha_l^b, v_1 \wedge \cdots \wedge v_l> = <\alpha_1 \wedge \cdots \wedge \alpha_l, v_1^{\#} \wedge \cdots \wedge v_l^{\#}> \\ (\alpha_1 \vee \cdots \vee \alpha_l)(v_1, \cdots, v_l) &= <\alpha_1^b \vee \cdots \vee \alpha_l^b, v_1 \vee \cdots \vee v_l> = <\alpha_1 \vee \cdots \vee \alpha_l, v_1^{\#} \vee \cdots \vee v_l^{\#}>\end{aligned} \tag{17}$$

Now, for $\xi = \alpha_1 \vee \cdots \vee \alpha_{l+1}$ and $\eta = \beta_1 \vee \cdots \vee \beta_l$ the following computation proves the second equality.

$$< i_v \xi, \eta > = < i_v(\alpha_1 \vee \cdots \vee \alpha_{l+1}), \beta_1 \vee \cdots \vee \beta_l > = i_v(\alpha_1 \vee \cdots \vee \alpha_{l+1})(\beta_1^b, \cdots, \beta_l^b) \\ = (\alpha_1 \vee \cdots \vee \alpha_{l+1})(v, \beta_1^b, \cdots, \beta_l^b) = < \xi, v^\# \vee \eta > \quad (18)$$

The same computation shows the first equality. ●

For any base $\{e_1, \cdots, e_n\}$ (also denoted by $\{e_i\}$) its dual base is denoted by $\{e^i\}$. The isomorphism $\alpha \mapsto \alpha^b$, (for symmetric or symplectic product) the dual base $\{e^i\}$ is transformed to the base $\{(e^i)^b\}$. This base is called the reciprocal base of $\{e_i\}$. The characteristic property of the reciprocal base is $< (e^j)^b, e_i > = \delta_i^j$. To prevent cumbersome notation, we denote $(e^i)^b$ by $e^i$, too. So $e^i$ denote a vector in $W$ or a vector in $W^*$, depending on the context. So, for any base $\{e_1, \cdots, e_n\}$, its reciprocal base or dual base is $\{e^1, \cdots, e^n\}$ and $< e^j, e_i > = \delta_i^j$, also $e^j(e_i) = \delta_i^j$.

**Super inner product:** A super inner product on a super space $V$, is a scalar valued nonsingular bilinear map on $V$, such that:

$$< u, v > = (-1)^{uv} < v, u > \\ u \in V_0, v \in V_1 \Rightarrow < u, v > = 0 \quad (19)$$

The second property means that inner product preserves parity, and the first property means that this product is an (ordinary) inner product on $V_0$ and is a nonsingular symplectic product on $V_1$.

A base of $V$ whose vectors are homogeneous is called a homogeneous base of $V$. A homogeneous base of $V$ is union of a base of $V_0$ and a base of $V_1$. Nondegeneracy of the super inner product, implies that for every homogeneous bases $\{e_i\}$ there exists a unique base $\{e^i\}$ such that $< e^j, e_i > = \delta_i^j$. The base $\{e^i\}$ is homogenous too, and $|e_i| = |e^i|$. The base $\{e^i\}$ is called the reciprocal base of $\{e_i\}$. The reciprocal base of $\{e^i\}$ is $\{(-1)^{e_i} e_i\}$. The dual base of $\{e_i\}$ is also denoted by $\{e^i\}$.

**Super contraction of tensors:** Suppose $V$ is a super space with a super inner product and $T: V \times V \to \mathbb{R}$ is a bilinear map that preserves parity. The super contraction of $T$ is a scalar, corresponds to the super trace of linear maps.

To define super contraction of $T$, consider a homogenous base $\{e_i\}$ and its reciprocal $\{e^i\}$. The super contraction of $T$ is define as follows.

$$\text{Super contraction of } T = \text{sc}(T) = \sum_i T(e_i, e^i) \tag{20}$$

This definition is well defined and does not depend on the choice of the base. If $T$ is super alternating (i.e. $(u,v) = -(-1)^{uv} T(v,u)$ ) then its super contraction is zero. So, the super contraction of a tensor $T$ only depends on its super symmetric part.

**Hodge operator:** In this section, we remind some basic properties of the Hodge operator in vector spaces equipped with an ordinary inner product[1]. Those vector $v$ that $<v,v>$ is negative (positive), are called negative (positive) vectors.

If $\{e_1, \cdots, e_n\}$ is an orthonormal base (i.e. $<e_i, e_j> = \pm \delta_{ij}$) then its dual base $\{e^1, \cdots, e^n\}$ is also orthonormal moreover, $e_i$ and $e^i$ have the same sign.

For $\alpha \in W^*$ the insertion operator $i_\alpha$ can be defined by $i_\alpha = i_{\alpha^b}$. This operation can be extended for $\eta \in \wedge^l W^*$. For a simple tensor $\eta = \beta_1 \wedge \cdots \wedge \beta_l$, $i_\eta(\omega)$ is defined by $i_\eta = i_{\beta_1 \wedge \cdots \wedge \beta_l} = i_{\beta_l} \circ \cdots \circ i_{\beta_1}$.

Assume $W$ is oriented. This orientation induce an orientation on $W^*$. If $\{e_1, \cdots, e_n\}$ is positive oriented then its dual, $\{e^1, \cdots, e^n\}$ is also positive oriented.

The canonical volume tensor $\Omega \in \wedge^n W^*$ is defined by $\Omega = e^1 \wedge \cdots \wedge e^n$ in which $\{e_1, \cdots, e_n\}$ is a positive oriented orthonormal base of $W$.

For every $k = 0, \cdots, n$, the Hodge operator $H: \wedge^k W^* \to \wedge^{n-k} W^*$ is defined by $H(\eta) = i_\eta(\Omega)$. It can be proved that for any transformation $\sigma \in S_n$:

$$H(e^{\sigma(1)} \wedge \cdots \wedge e^{\sigma(k)}) = (-1)^t \varepsilon_\sigma e^{\sigma(k+1)} \wedge \cdots \wedge e^{\sigma(n)} \tag{21}$$

In above, $t$ is the number of negative vectors in $e^{\sigma(1)}, \cdots, e^{\sigma(k)}$. If the number of negative vector in $\{e_1, \cdots, e_n\}$ is $s$, then

$$H \circ H = (-1)^{k(n-k)+s} \mathbf{1}, \quad H(1) = \Omega, \ H(\Omega) = (-1)^s \tag{22}$$

For $\omega \in \wedge^k W^*$ and $\eta \in \wedge^l W^*$ we have $H(\omega \wedge \eta) = i_\eta(\omega)$.

An important property of the Hodge operator is that for all $\omega, \varphi \in \Lambda^k W^*$ the following equality holds.

$$\omega \wedge H(\varphi) = <\omega, \varphi> \Omega \tag{23}$$

Now, assume $V$ is a super space that is equipped with a super inner product. We extend Hodge operator on super alternating tensors. First, we consider the ordinary Hodge operator on $\Lambda^k V_0^*$. For simple tensors such as $\omega \otimes \xi$, it is defined:

$$H(\omega \otimes \xi) = H(\omega) \otimes \xi \tag{24}$$

The symmetric product on $V_0$ and the symplectic product on $V_1$, induce some products on $\Lambda^s V_0^*$,s and $\vee^l V_1^*$,s. These products, also induce some products on $(\Lambda^s V_0^*) \otimes (\vee^l V_1^*)$. For simple tensors $\omega \otimes \xi$ and $\varphi \otimes \eta$, their production is defined as follows.

$$<\omega \otimes \xi, \varphi \otimes \eta> = <\omega, \varphi><\xi, \eta> \tag{25}$$

Associated to the (symmetric or symplectic) product on $\vee^l V_1^*$, we can define an exterior product on $\hat{A}(V)$ with value in $A(V_0)$ that is denoted by $\dot\wedge$. For simple tensor $\omega \otimes \xi$ and $\varphi \otimes \eta$, it is defined:

$$(\omega \otimes \xi) \dot\wedge (\varphi \otimes \eta) = (\omega \wedge \varphi) <\xi, \eta> \tag{26}$$

By this definition, for $\tilde{\omega}, \tilde{\varphi} \in \hat{A}(V)$ we find that:

$$\tilde{\omega} \dot\wedge H(\tilde{\varphi}) = <\tilde{\omega}, \tilde{\varphi}> \Omega \tag{27}$$

These concepts can also be extended to vector valued alternating tensors. Suppose that $W$ is a super space with some super inner product on it. We can extend these inner products on $\hat{A}(V, W) = \hat{A}(V) \otimes W$. Hodge operator can be extended to $\hat{A}(V, W)$ by relation $H(\tilde{\omega} \otimes u) = H(\tilde{\omega}) \otimes u$., The exterior product on $\hat{A}(V, W)$ with value in $A(V_0)$, associated to the inner product on $W$, is denoted by $\dot\wedge$.

$$(\tilde{\omega} \otimes u) \dot\wedge (\tilde{\varphi} \otimes v) = (\tilde{\omega} \dot\wedge \tilde{\varphi}) <u, v> \tag{28}$$

These definitions imply that for $\boldsymbol{\omega}, \boldsymbol{\varphi} \in \hat{A}(V, W)$ we have:

$$\boldsymbol{\omega} \dot\wedge H(\boldsymbol{\varphi}) = <\boldsymbol{\omega}, \boldsymbol{\varphi}> \Omega \tag{29}$$

## 3. Super tangent bundles

Here, $M$ is a smooth manifold and all vector bundles are on $M$. A vector bundle $E$ with a given decomposition $E = E_0 \oplus E_1$ is called a super vector bundle. To make a super tangent bundle for $M$, we can add to $TM$ a vector bundle $F$. So, for a vector bundle $F$ set $\hat{T}M = TM \oplus F$. Clearly, $\hat{T}M$ depends on the choice of $F$. Even part of $\hat{T}M$ is $TM$ and the odd part is $F$. The set of all sections of $\hat{T}M$ is denoted by $\widehat{\mathfrak{X}}M$, that is equal to $\mathfrak{X}M \oplus \Gamma F$. Ordinary vector fields on $M$ are even vector fields and sections of $F$ are odd vector fields.

From now on, we use symbols $U, V, \cdots$ for ordinary vector fields, and $\alpha, \beta, \cdots$ for odd vector fields, and $\boldsymbol{U}, \boldsymbol{V}, \ldots$ for arbitrary vector fields. Cleary, every arbitrary vector field $\boldsymbol{U}$ can be written in the form $\boldsymbol{U} = U + \alpha$.

Action of odd vector fields on smooth functions is defined to be zero, because we have no odd function except zero. So, for $\alpha \in \Gamma F$ and $f \in C^\infty(M)$ we have $\alpha(f) = 0$. The most important property of tangent bundles is Lie bracket of vector fields. So, to have a super tangent bundle, it must be given a super Lie algebra structure on $\widehat{\mathfrak{X}}M$ such that its even Lie algebra is $\mathfrak{X}M$.

**Definition(3.1):** For a manifold $M$ and a vector bundle $F$ on $M$, the super vector bundle $\hat{T}M = TM \oplus F$ with a super Lie algebra structure on $\widehat{\mathfrak{X}}M$ is called a super tangent bundle on $M$, if its even Lie algebra is $\mathfrak{X}M$ with usual Lie bracket, and for $\boldsymbol{U}, \boldsymbol{V} \in \widehat{\mathfrak{X}}M$ and $f \in C^\infty(M)$ the following relation holds.

$$[\boldsymbol{U}, f\boldsymbol{V}] = \boldsymbol{U}(f)\boldsymbol{V} + f[\boldsymbol{U}, \boldsymbol{V}] \tag{30}$$

**Theorem(3.2):** Lie bracket of all odd vector fields in any super tangent bundle is zero.

**Proof:** It suffices to show that $[\alpha, \alpha] = 0$ for all odd vector field $\alpha$. We know that $[\alpha, [\alpha, \alpha]] = 0$. Also for every smooth function $f$, $[f\alpha, [f\alpha, f\alpha]] = 0$. So,

$$\begin{aligned} 0 = [f\alpha, [f\alpha, f\alpha]] &= [f\alpha, f^2[\alpha, \alpha]] \\ &= f^2(-[\alpha, \alpha](f)\alpha + f[\alpha, [\alpha, \alpha]]) = -f^2[\alpha, \alpha](f)\alpha \end{aligned} \tag{31}$$

Since $f$ is arbitrary, we can deduce $[\alpha, \alpha](f) = 0$ and consequently, $[\alpha, \alpha] = 0 \bullet$

**Theorem(3.3):** If $M$ is simply connected, and $\hat{T}M = TM \oplus F$ is a super tangent bundle, then $F$ is a trivial bundle and the Lie bracket of an ordinary vector field and an odd vector field is the action of a vector fields on a vector valued function.

**Proof:** For $U \in \mathfrak{X}M$ and $\alpha \in \Gamma F$ define: $\nabla_U \alpha = [U, \alpha]$. The properties of Lie bracket imply that $\nabla$ is a connection on the bundle $F$. From Jacoby identity for two ordinary and one odd vector field, we find the curvature tensor of this connection is zero. Because,

$$R(U,V)(\alpha) = \nabla_U \nabla_V \alpha - \nabla_V \nabla_U \alpha - \nabla_{[U,V]} \alpha \\ = [U,[V,\alpha]] - [V,[U,\alpha]] - [[U,V],\alpha] = 0 \tag{32}$$

Since $M$ is simply connected, parallel translation of odd vectors between two points is independent of the chosen path[6]. So, by parallel translation, we can construct a trivialization for $F$, and with respect to this trivialization, $\nabla$ is trivial. In other words, we can assume, $F = M \times W$ for some vector spaces $W$ and we have

$$[U, \alpha] = U(\alpha). \; \bullet \tag{33}$$

Due to this theorem, every super tangent bundle on a simply connected manifold, is determined by a vector space $W$ and $\hat{T}M = TM \oplus (M \times W)$ and odd fields are $W$-valued smooth functions on $M$, and for odd fields $\alpha, \beta$ and ordinary vector fields $U, V$, we have

$$[U + \alpha, V + \beta] = [U, V] + U(\beta) - V(\alpha) \tag{34}$$

We investigate only this case, because most of manifolds are simply connected and all manifolds are locally simply connected. In this case, some odd vector fields are constant function and are called constant odd vector fields. Clearly, $\alpha$ is a constant odd vector field if and only if for every ordinary vector field $U$, $[U, \alpha] = 0$. For odd fields, we can choose a base of constant odd vector fields.

4. **super differential forms**

From now on, we consider a fixed super tangent bundle $\hat{T}M = TM \oplus (M \times W)$. For a super bundle $E$, a covariant $E$-valued (super) tensor field on $M$ can be defined as an operator that is $C^\infty(M)$-multilinear.

$$\omega: \widehat{\mathfrak{X}}M \times \cdots \times \widehat{\mathfrak{X}}M \to \Gamma E \tag{35}$$

If $\omega$ preserve parities of vectors, it is called even, and if it reverse parities of vectors, it is called odd. For special case $= M \times \mathbb{R}$, $E$-valued tensors are called numerical valued tensors.

The tensor field $\omega$ of order $k$, is called (super) alternating, if

$$\omega(U_1, \cdots, U_i, U_{i+1}, \cdots, U_k) = -(-1)^{U_i U_{i+1}} \omega(U_1, \cdots, U_{i+1}, U_i, \cdots, U_k) \tag{36}$$

These tensors are called (super) differential forms. The set of all $E$-valued (super) $k$-differential forms is denoted by $\hat{A}^k(M, E)$. For the special case of numerical valued $k$-(super) differential forms, this set is denoted by $\hat{A}^k(M)$.

For $\omega \in A^s(M)$ and $\xi \in C^\infty(M, \vee^l W^*)$, we can define $\omega \otimes \xi$ as a numerical super alternating tensor of order $s + l$ as follows.

$$(\omega \otimes \xi)(U_1, \cdots, U_s, \alpha_1, \cdots, \alpha_l) = \omega(U_1, \cdots, U_s)\xi(\alpha_1, \cdots, \alpha_l) \tag{37}$$

We call these tensors, simple tensors and the pair $(s, l)$ is called the pair order of $\omega \otimes \xi$. So, (super) differential forms have two parities. The second parity of these tensors is identical to be even or odd, as defined above. We use the same notation for the first and second parity of tensor fields. For a simple super differential form $\omega \otimes \xi$, $\omega$ is called its alternating part and $\xi$ is called its symmetric part.

Super differential forms are sum of simple tensors; moreover we can choose those simple tensors whose symmetric parts are constant.

Exterior product of (super) differential forms is defined similarly as ones on the alternating tensors on super spaces. Indeed, for simple differential forms, we have

$$(\omega \otimes \xi) \wedge (\varphi \otimes \eta) = (\omega \wedge \varphi) \otimes (\xi \vee \eta) \tag{38}$$

This product preserves first and second parities.

For a super bundle $E$ and $X \in \Gamma E$ and $\tilde{\omega} \in \hat{A}^k(M)$, we define $\tilde{\omega} \otimes X$ to be the following $E$-valued differential form.

$$(\tilde{\omega} \otimes X)(U_1, \cdots, U_k) = \tilde{\omega}(U_1, \cdots, U_k)X \tag{39}$$

Associated to any product (bilinear map preserving parity) between super vector bundles it is defined an exterior product between vector valued differential forms. Exterior product of simple vector valued differential forms is as follows.

$$(\widetilde{\omega} \otimes X) \wedge (\widetilde{\eta} \otimes Y) = (\widetilde{\omega} \wedge \widetilde{\eta}) \otimes (XY) \tag{40}$$

First and second parity of vector valued tensor fields are defined similarly and exterior products preserve both parities.

## 5. Basic differential operators

In this section, we explain some basic derivations on the super algebra $\hat{A}(M)$ (with respect to the first parity)[4][5]. A homogenous linear map $D$ on $\hat{A}(M)$ is called a (super) derivation if

$$D(\widetilde{\omega} \wedge \widetilde{\varphi}) = (D\widetilde{\omega}) \wedge \widetilde{\varphi} + (-1)^{\widetilde{\omega}D} \widetilde{\omega} \wedge D\widetilde{\varphi} \tag{41}$$

If $D_1$ and $D_2$ are two homogeneous (super) derivation then $[D_1, D_2]$ that is defined as follows, is a derivation whose parity is $|D_1| + |D_2|$.

$$[D_1, D_2] = D_1 \circ D_2 - (-1)^{D_1 D_2} D_2 \circ D_1 \tag{42}$$

**(Super) Lie derivations**: For $U \in \widehat{\mathfrak{X}}M$ it is constructed the operator $L_U$ on $\hat{A}(M)$ as follows. For example, for a super differential form $\widetilde{\omega}$ of order two, we have:

$$(L_U \widetilde{\omega})(V, W) = U(\widetilde{\omega}(V, W)) - (-1)^{\widetilde{\omega} U} \widetilde{\omega}([U, V], W) - (-1)^{(\widetilde{\omega}+V)U} \widetilde{\omega}(V, [U, W]) \tag{43}$$

This definition is well defined. The operator $L_U$ is called (super) Lie derivation along $U$ and does not change the order of tensors. If we consider the first parity for (super) differential forms, then the parity of $L_U$ is the same as the parity of $U$. In other words, for even vector fields $U$, the operator $L_U$ does not change the first parity, and for odd vector fields $\alpha$, the operator $L_\alpha$ reverse the first parity. For simple differential forms $\omega \in A^k(M)$ and $\xi \in C^\infty(M, \vee^l W^*)$ and $L_U \omega$ and $L_U \xi$ are ordinary Lie derivation, but $L_\alpha \omega = 0$ and the pair order of $L_\alpha \xi$ is $(1, l-1)$.

$$(L_\alpha \xi)(U, \beta_1, \cdots, \beta_{l-1}) = \xi(U(\alpha), \beta_1, \cdots, \beta_{l-1}) \tag{44}$$

If $\alpha$ is constant, then $L_\alpha \xi = 0$. A straightforward computation shows that

$$L_U(\omega \otimes \xi) = (L_U\omega) \otimes \xi + (-1)^{\omega U} \omega \wedge L_U \xi \tag{45}$$

This relation shows that in general case for two (super) differential forms $\tilde{\omega}$ and $\tilde{\varphi}$ the following equality holds.

$$L_U(\tilde{\omega} \wedge \tilde{\varphi}) = (L_U\tilde{\omega}) \wedge \tilde{\varphi} + (-1)^{\tilde{\omega}U} \tilde{\omega} \wedge L_U\tilde{\varphi} \tag{46}$$

So, $L_U$ is a super derivation on $\hat{A}(M)$.

**Insertion operators:** For $U \in \widehat{\mathfrak{X}}M$ the operator $i_U$ on $\hat{A}(M)$ is an operator that reduce order of tensors. For $\tilde{\omega} \in \hat{A}^k(M)$, the order of $i_U\tilde{\omega}$ is $k-1$ and defined as follows.

$$(i_U\tilde{\omega})(V_1, \cdots, V_{k-1}) = (-1)^{\tilde{\omega}U} \tilde{\omega}(U, V_1, \cdots, V_{k-1}) \tag{47}$$

For smooth functions (zero order tensors) such as $f$, it is defined $i_U(f) = 0$. If $U$ is homogeneous, then the operator $i_U$ is also homogenous and with respect to the first parity of differential forms, its parity is $|U| + 1$. In other words, the operator $i_U$ reverse the first parity and the operator $i_\alpha$ preserve the first parity. For simple tensors such as $\omega \otimes \xi$ we have:

$$i_U(\omega \otimes \xi) = (i_U\omega) \otimes \xi \quad , \quad i_\alpha(\omega \otimes \xi) = \omega \otimes i_\alpha \xi \tag{48}$$

This relation shows that for two (super) differential forms $\tilde{\omega}$ and $\tilde{\varphi}$ we have:

$$i_U(\tilde{\omega} \wedge \tilde{\varphi}) = (i_U\tilde{\omega}) \wedge \tilde{\varphi} + (-1)^{\tilde{\omega}(U+1)} \tilde{\omega} \wedge i_U\tilde{\varphi} \tag{49}$$

So, $i_U$ is a super derivation on $\hat{A}(M)$.

**(Super) exterior derivation:** Naturally, exterior derivation of ordinary differential forms can be extended to a super derivation on $\hat{A}(M)$ that we call it, super exterior derivation. This operator is denoted by $d$ and increases the order of tensors. For $\tilde{\omega} \in \hat{A}^k(M)$, the order of $d\tilde{\omega}$ is $k+1$ and defined as follows.

$$\begin{aligned}(d\tilde{\omega})(U_1, \cdots, U_{k+1}) &= \sum_{i=1}^{k+1}(-1)^{i+1+U_i(U_1+\cdots+U_{i-1})} U_i(\tilde{\omega}(U_1, \cdots, \hat{U}_i, \cdots, U_{k+1})) \\ &+ \sum_{1 \le i < j \le k+1}(-1)^{i+j+U_i(U_1+\cdots+U_{i-1})+U_j(U_1+\cdots+\hat{U}_i+\cdots+U_{j-1})} \tilde{\omega}([U_i, U_j], U_1, \cdots, \hat{U}_i, \cdots, \hat{U}_j, \cdots, U_{k+1})\end{aligned} \tag{50}$$

By a straightforward computation we find that $d\tilde{\omega}$ is a tensor. So, in computations, we can assume the Lie brackets of vector fields in above are zero and prove that $d\widetilde{\omega}$ is (super) alternating. This operator reverses the first parity of tensors and is an odd operator.

An ordinary differential form $\omega \in A^k(M)$ as a super differential form has the pair order $(k, 0)$ and $d\omega$ is its ordinary exterior derivation whose pair order is $(k+1, 0)$.

But, a tensor $\xi \in C^\infty(M, \vee^l W^*)$ as a super differential form has the pair order $(0, l)$ and $d\xi$ has the pair order $(1, l)$. We can easily show that:

$$d\xi(U, \beta_1, \cdots, \beta_l) = U(\xi)(\beta_1, \cdots, \beta_l) \tag{51}$$

If $\xi$ is constant, then $d\xi = 0$. For a simple tensor $\omega \otimes \xi$ we find that:

$$d(\omega \otimes \xi) = d\omega \otimes \xi + (-1)^\omega \omega \wedge d\xi \tag{52}$$

Especially, if $\xi$ is constant we have $d(\omega \otimes \xi) = d\omega \otimes \xi$. By this relation we can deduce that $d \circ d = 0$.

In general, if the pair order of $\widetilde{\omega}$ is $(s, l)$ then the pair order of $d\widetilde{\omega}$ is $(s+1, l)$. Super exterior derivation is a super derivation and the following relation holds.

$$d(\widetilde{\omega} \wedge \widetilde{\varphi}) = (d\widetilde{\omega}) \wedge \widetilde{\varphi} + (-1)^{\widetilde{\omega}} \widetilde{\omega} \wedge d\widetilde{\varphi} \tag{53}$$

We can compute Lie brackets of these three super derivations.

For two insertion operator $i_U$ and $i_V$ we find that $[i_U, i_V] = 0$. For example, for a differential form $\widetilde{\omega}$ of order two we have:

$$\begin{aligned}
[i_U, i_V](\widetilde{\omega}) &= \left(i_U i_V - (-1)^{(U+1)(V+1)} i_V i_U\right)(\widetilde{\omega}) \\
&= (-1)^{U(\widetilde{\omega}+V+1)} i_V(\widetilde{\omega})(U) - (-1)^{(U+1)(V+1)+V(\widetilde{\omega}+U+1)} i_U(\widetilde{\omega})(V) \\
&= (-1)^{U(\widetilde{\omega}+V+1)+V\widetilde{\omega}} \widetilde{\omega}(V, U) - (-1)^{U+1+V\widetilde{\omega}+U\widetilde{\omega}} \widetilde{\omega}(U, V) \\
&= (-1)^{U+U\widetilde{\omega}+V\widetilde{\omega}} \left((-1)^{UV} \widetilde{\omega}(V, U) + \widetilde{\omega}(U, V)\right) = 0
\end{aligned} \tag{54}$$

Lie bracket of the exterior derivation and insertion operators are Lie derivations. In other words,

$$[d, i_U] = d \circ i_U - (-1)^{U+1} i_U \circ d = L_U \tag{55}$$

To prove this equality, we can prove it for even and odd vector fields separately. And in each case, it is better to prove it first, for ordinary tensors $\omega \in A^k(M)$ and $\xi \in C^\infty(M, \vee^l W^*)$ and then for simple tensors $\omega \otimes \xi$.

By these relations we find that $[d, L_U] = 0$. By straightforward computations, we also find that $[L_U, L_V] = L_{[U,V]}$.

## 6. Super connections on super vector bundles

Let $E$ be super vector bundle on $M$, and $\widehat{T}M$ is a super tangent bundle on $M$.

**Definition(6.1):** A super connection (or super covariant derivation) on $E$ is a bilinear map (with respect to scalars)

$$\nabla: \widehat{\mathfrak{X}}M \times \Gamma E \longrightarrow \Gamma E, \quad (U, X) \longmapsto \nabla_U X \tag{56}$$

such that the following relations hold.

1) $|\nabla_U X| = |U| + |X|$ \hfill (57)
2) $\nabla_U(fX) = U(f)X + f\nabla_U X$ \hfill (58)
3) $\nabla_{fU} X = f\nabla_U X$ \hfill (59)

The restriction of $\nabla$ to $\mathfrak{X}M \times \Gamma E$ is an ordinary connection on $E$, that preserves even and odd sub bundles. But, the restriction of $\nabla$ to odd fields is a tensor $\mathcal{D}(\alpha, X)$ that with respect to $X$ is an odd bundle map. Therefore, a super connection on $E$ is determined by a pair $(\nabla, \mathcal{D})$ in which, $\nabla$ is an ordinary connection on $E$, that preserves even and odd sub bundles and $\mathcal{D}$ is a tensor with the properties described above.

**Example(6.2):** For the case of trivial super bundle, the operator $(U, X) \longmapsto U(X)$ is a super connection. In this case, $\mathcal{D} = 0$.

In the case $E = \widehat{T}M$, $\nabla$ is called a super connection on $M$. The restriction of $\nabla$ to $\mathfrak{X}M \times \mathfrak{X}M$ is an ordinary connection on $M$, and its restriction to $\mathfrak{X}M \times C^\infty(M, W)$ is a connection on the trivial bundle $M \times W$. So, for some $\omega_0 \in A^1(M, L(W))$, we have:

$$\nabla_U \alpha = U(\alpha) + \omega_0(U, \alpha) \tag{60}$$

Moreover, for every odd field $\alpha$ the operator $\omega_\alpha: U \mapsto \nabla_\alpha U$ is a $W$-valued one form on $M$. In other words, $\omega_\alpha \in A^1(M, W)$, furthermore, $\alpha \mapsto \omega_\alpha$ is linear. For any two odd fields $\alpha, \beta$, $\nabla_\alpha \beta$ is an ordinary vector field that we denote it by $X_{\alpha\beta}$. Clearly, $(\alpha, \beta) \mapsto X_{\alpha\beta}$ is bilinear. A super connection on $M$ is determined by an ordinary connection on $M$, and some tensors $\omega_0, \omega_\alpha$, and $X_{\alpha\beta}$.

The (super)torsion of a super connection on $M$, is a $\hat{T}M$-valued (super) differential form of order two, and is defined as follows.

$$T(U, V) = \nabla_U V - (-1)^{UV} \nabla_V U - [U, V] \tag{61}$$

Since, $T$ preserve parities of vector fields($|T(U,V)| = |U| + |V|$), it is an even tensor. The super torsion of a super connection is zero, if and only if its associated ordinary connection is torsion free and

$$\omega_0(U, \alpha) = \omega_\alpha(U), \quad X_{\alpha\beta} = -X_{\alpha\beta} \tag{62}$$

Therefore, a torsion free super connection $\nabla$ on $M$ is determined by a triplet $\nabla = (\dot{\nabla}, \omega_0, X_{\alpha\beta})$ in which $\dot{\nabla}$ a torsion free ordinary connection on $M$, and $\omega_0$ is a tensor in $A^1(M, L(W))$, and $X_{\alpha\beta}$ are vector fields such that with respect to $\alpha, \beta$ is $C^\infty(M)$-bilinear and alternating. In this case, the covariant derivations between odd and even vector fields are as follow.

$$\nabla_U \alpha = U(\alpha) + \omega_0(U, \alpha), \quad \nabla_\alpha U = \omega_0(U, \alpha), \quad \nabla_\alpha \beta = X_{\alpha\beta} \tag{63}$$

Let $\nabla$ be a super connection on $M$. This connection can be extended to the tonsorial bundles associated to $\hat{T}M$. For example, if $E$ is another super vector bundle with some super connection on it, and $T$ is an $E$-valued super differential form on $M$ of order three, then $\nabla_U T$ is a tensor of the same type, and is defined by the following relation.

$$\nabla_U T(X, Y, Z) = (\nabla_U T)(X, Y, Z) + T(\nabla_U X, Y, Z) \\ + (-1)^{U.X} T(X, \nabla_U Y, Z) + (-1)^{U.(X+Y)} T(X, Y, \nabla_U Z) \tag{64}$$

The operator $\nabla_U$ does not changes order of tensors, moreover for an even vector field $U$, the operator $\nabla_U$ does not changes pair order of tensors. But, for an odd vector field $\alpha$ and a tensor $\tilde{\omega}$ of the pair order $(s,l)$, $\nabla_\alpha \tilde{\omega}$ is sum of some tensors of the pair orders $(s-1, l+1)$ and $(s+1, l-1)$.

For example, if $\omega \in A^k(M)$, then the pair order of $\omega$ is $(k, 0)$, and the pair order of $\nabla_\alpha \omega$ is $(k-1, 1)$, and we have:

$$(\nabla_\alpha \omega)(V_1, \cdots, V_{k-1}, \beta) = \omega(\nabla_\alpha \beta, V_1, \cdots, V_{k-1}) \tag{65}$$

Also, if $\xi \in C^\infty(M, V^l W^*)$, then the pair order of $\xi$ is $(0, l)$, and the pair order of $\nabla_\alpha \xi$ is $(1, l-1)$, and we have:

$$(\nabla_\alpha \xi)(V, \beta_1 \cdots, \beta_{l-1}) = -\xi(\nabla_\alpha V, \beta_1, \cdots, \beta_{l-1}) \tag{66}$$

In general, for a simple tensor $\omega \otimes \xi$ whose pair order is $(k, l)$, for an even vector field $U$ we have:

$$\nabla_U(\omega \otimes \xi) = (\nabla_U \omega) \otimes \xi + \omega \otimes \nabla_U \xi \tag{67}$$

Clearly, the pair order of $\nabla_U(\omega \otimes \xi)$ is also $(s, l)$. But, for an odd field $\alpha$, we have:

$$\nabla_\alpha(\omega \otimes \xi) = (\nabla_\alpha \omega) \wedge \xi + (-1)^s \omega \wedge \nabla_\alpha \xi \tag{68}$$

The pair order of $(\nabla_\alpha \omega) \wedge \xi$ is $(s-1, l+1)$, and the pair order of $\omega \wedge \nabla_\alpha \xi$ is $(s+1, l-1)$.

For local computations, we need some local bases of homogeneous vector fields. For odd vector fields, we can construct and fix a global base $\{\alpha_1, \cdots, \alpha_m\}$ of odd vector fields. Moreover, we can choose these vectors to be constant. The dual of this base is denoted by $\{\alpha^1, \cdots, \alpha^m\}$ and is constant provided that $\{\alpha_1, \cdots, \alpha_m\}$ is constant. For even vector fields, we consider some fixed local base $\{E_1, \cdots, E_n\}$ with dual base of 1-forms $\{E^1, \cdots, E^n\}$. Clearly, $\{E_1, \cdots, E_n, \alpha_1, \cdots, \alpha_m\}$ is a homogenous local base and $\{E^1, \cdots, E^n, \alpha^1, \cdots, \alpha^m\}$ is its dual base.

**Theorem(6.3)**: If $\omega \in A^k(M)$ then

$$\sum_j \alpha^j \wedge \nabla_{\alpha_j} \omega = 0 \tag{69}$$

**Proof**: Let $\tilde{\eta} = \sum_j \alpha^j \wedge \nabla_{\alpha_j} \omega$. The pair order of $\tilde{\eta}$ is $(k-1, 2)$. For $k = 0$, clearly $\tilde{\eta} = 0$. For $1 \leq k$ the following computations prove the equality.

$$\begin{aligned}
\tilde{\eta}(V_1,\cdots,V_{k-1},\beta,\gamma) &= \left(\sum_j \alpha^j \wedge \nabla_{\alpha_j}\omega\right)(V_1,\cdots,V_{k-1},\beta,\gamma) \\
&= \sum_j \alpha^j(\beta)\omega\left(\nabla_{\alpha_j}\gamma,V_1,\cdots,V_{k-1}\right) + \alpha^j(\gamma)\omega\left(\nabla_{\alpha_j}\beta,V_1,\cdots,V_{k-1}\right) \\
&= \omega(\nabla_\beta\gamma,V_1,\cdots,V_{k-1}) + \omega(\nabla_\gamma\beta,V_1,\cdots,V_{k-1}) = 0 \quad \bullet
\end{aligned} \tag{70}$$

Super exterior derivation of a differential form can be computed via a torsion free connection. In fact, the famous relation between ordinary exterior derivation and an ordinary torsion free connection also holds for super exterior derivation and a torsion free super connection.

**Theorem(6.4):** If $\tilde{\omega} \in \hat{A}(M)$, then $d\tilde{\omega} = \sum E^i \wedge \nabla_{E_i}\tilde{\omega} - \sum \alpha^j \wedge \nabla_{\alpha_j}\tilde{\omega}$.

**Proof:** First, we prove this equality for $\omega \in A^k(M)$ and $\xi \in C^\infty(M, \vee^l W^*)$. Since $\sum \alpha^i \wedge \nabla_{\alpha_i}\omega = 0$, the equality holds for $\omega$. Some lengthy computations prove the equality for $\xi$.

$$\begin{aligned}
\left(\sum_i E^i \wedge \nabla_{E_i}\xi - \sum_j \alpha^j \wedge \nabla_{\alpha_j}\xi\right)(V,\beta_1,\cdots,\beta_l) & \\
&= \sum_i E^i(V)(\nabla_{E_i}\xi)(\beta_1,\cdots,\beta_l) - \sum_{j,t} \alpha^j(\beta_t)\left(\nabla_{\alpha_j}\xi\right)(V,\beta_1,\cdots,\hat{\beta}_t,\cdots,\beta_l) \\
&= (\nabla_V\xi)(\beta_1,\cdots,\beta_l) - \sum_t (\nabla_{\beta_t}\xi)(V,\beta_1,\cdots,\hat{\beta}_t,\cdots,\beta_l) \\
&= V(\xi(\beta_1,\cdots,\beta_l)) - \sum_t \xi(\beta_1,\cdots,\nabla_V\beta_t,\cdots,\beta_l) + \sum_t \xi(\nabla_{\beta_t}V,\beta_1,\cdots,\hat{\beta}_t,\cdots,\beta_l) \\
&= V(\xi)(\beta_1,\cdots,\beta_l) + \sum_t \xi(\beta_1,\cdots,V(\beta_t),\cdots,\beta_l) + \sum_t \xi(\beta_1,\cdots,\nabla_{\beta_t}V - \nabla_V\beta_t,\cdots,\beta_l) \\
&= d\xi(V,\beta_1,\cdots,\beta_l) + \sum_t \xi(\beta_1,\cdots,V(\beta_t),\cdots,\beta_l) + \sum_t \xi(\beta_1,\cdots,-V(\beta_t),\cdots,\beta_l) \\
&= d\xi(V,\beta_1,\cdots,\beta_l)
\end{aligned} \tag{71}$$

Now, we are ready to prove the equality for a simple tensor $\omega \otimes \xi$.

$$\begin{aligned}
\sum_i E^i \wedge \nabla_{E_i}(\omega \otimes \xi) &- \sum_j \alpha^j \wedge \nabla_{\alpha_j}(\omega \otimes \xi) \\
&= \sum_i E^i \wedge \left((\nabla_{E_i}\omega) \otimes \xi + \omega \otimes \nabla_{E_i}\xi\right) - \sum_j \alpha^j \wedge \left(\left(\nabla_{\alpha_j}\omega\right) \wedge \xi + (-1)^k \omega \wedge \nabla_{\alpha_j}\xi\right) \\
&= \sum_i \left(E^i \wedge (\nabla_{E_i}\omega)\right) \otimes \xi + (-1)^k \omega \wedge \sum_i E^i \wedge \nabla_{E_i}\xi \\
&\quad - \sum_j (\alpha^j \wedge (\nabla_{\alpha_j}\omega)) \wedge \xi - (-1)^k \omega \wedge \left(\sum_j \alpha^j \wedge \nabla_{\alpha_j}\xi\right) \\
&= d\omega \otimes \xi + (-1)^k \omega \wedge \left(\sum_i E^i \wedge \nabla_{E_i}\xi - \sum_j \alpha^j \wedge \nabla_{\alpha_j}\xi\right) \\
&= d\omega \otimes \xi + (-1)^\omega \omega \wedge d\xi = d(\omega \otimes \xi) \quad \bullet
\end{aligned} \tag{72}$$

**Super curvature:** For a super connection on a super vector bundle, its super curvature naturally is defined as follows.

$$R(U,V)(X) = \nabla_U\nabla_V X - (-1)^{UV}\nabla_V\nabla_U X - \nabla_{[U,V]}X \tag{73}$$

The super curvature is a super tensor that preserves the parities and is super alternating with respect to $U$ and $V$. In other words $R$ is a $L(E)$-valued super differential form of order two.

$$R(U,V)(X) = -(-1)^{UV}R(V,U)(X) \tag{74}$$

In the case $E = \hat{T}M$ for torsion free connections, Bianchi identities hold as follows.

$$\begin{aligned}(-1)^{UW}R(U,V)(W) + (-1)^{VU}R(V,W)(U) + (-1)^{WV}R(W,U)(V) &= 0 \\ (-1)^{UW}(\nabla_U R)(V,W) + (-1)^{VU}(\nabla_V R)(W,U) + (-1)^{WV}(\nabla_W R)(U,V) &= 0\end{aligned} \tag{75}$$

Proof is the same as ordinary Bianchi identities and in computations we can assume Lie brackets of vector fields (even or odd) are zero.

### 7. Basic differential operators associated to a super connection

In this section, we extend differential operators defined in section (5) to vector valued differential forms[5]. Assume $E$ is a super vector bundle equipped with a super connection $\nabla$. For a $U \in \hat{\mathfrak{X}}M$ we can extend Lie derivation $L_U$ on $\hat{A}(M)$ to $L_U^\nabla: \hat{A}(M,E) \to \hat{A}(M,E)$ as follows. For simple tensors $\tilde{\omega}\otimes X$, set:

$$L_U^\nabla(\tilde{\omega}\otimes X) = (L_U\tilde{\omega})\otimes X + (-1)^{\tilde{\omega}U}\tilde{\omega}\otimes\nabla_U X \tag{76}$$

With respect to the first parity, if $U$ is homogeneous, then the operator $L_U^\nabla$ is homogeneous and its parity is the same as the parity of $U$. If it is given a product on $E$ that is parallel with respect to the connection, then with respect to the associated exterior product, the operator $L_U^\nabla$ is a super derivation. In other words for $\omega, \varphi \in \hat{A}(M,E)$ we have:

$$L_U^\nabla(\omega \wedge \varphi) = (L_U^\nabla\omega) \wedge \varphi + (-1)^{\omega U}\omega \wedge L_U^\nabla\varphi \tag{77}$$

For different super vector bundles and different super connections on them and a product between them that is parallel with respect to that super connections, the same relation holds for the associated exterior product. For example, if we consider the trivial bundle $M \times \mathbb{R}$ and the trivial connection on it, then scalar product between $M \times \mathbb{R}$ and super bundle $E$, is parallel. So for associated exterior product, and $\tilde{\omega} \in \hat{A}(M)$ and $\varphi \in \hat{A}(M, E)$ we have:

$$L_U^\nabla(\tilde{\omega} \wedge \varphi) = (L_U\tilde{\omega}) \wedge \varphi + (-1)^{\tilde{\omega}U}\tilde{\omega} \wedge L_U^\nabla\varphi \tag{78}$$

Extension of insertion operators needs no super connection on $E$. For a simple tensor $\tilde{\omega} \otimes X$ it is defined:

$$i_U(\tilde{\omega} \otimes X) = (i_U\tilde{\omega}) \otimes X \tag{79}$$

This extended insertion operators have the same properties as the ordinary insertion operators. Especially, with respect of an exterior product between super vector bundles, we have:

$$i_U(\omega \wedge \varphi) = (i_U\omega) \wedge \varphi + (-1)^{\omega(U+1)}\omega \wedge i_U\varphi \tag{80}$$

Associated to the super connection on $E$, it is defined the exterior derivation on $E$-valued differential forms $d^\nabla: \hat{A}^k(M, E) \to \hat{A}^{k+1}(M, E)$ as follows.

$$\begin{aligned}(d^\nabla\omega)(U_1, \cdots, U_{k+1}) &= \sum_{i=1}^{k+1}(-1)^{i+1+U_i(U_1+\cdots+U_{i-1})}\nabla_{U_i}\omega(U_1, \cdots, \hat{U}_i, \cdots, U_{k+1}) \\ &+ \sum_{1 \le i < j \le k+1}(-1)^{i+j+U_i(U_1+\cdots+U_{i-1})+U_j(U_1+\cdots+\hat{U}_i+\cdots+U_{j-1})}\omega([U_i, U_j], U_1, \cdots, \hat{U}_i, \cdots, \hat{U}_j, \cdots, U_{k+1})\end{aligned} \tag{81}$$

This definition is well-defined and the proof is similar to the case of ordinary exterior derivation. For example, if $\omega \in \hat{A}^1(M, E)$, then

$$(d^\nabla\omega)(U, V) = \nabla_U\omega(V) - (-1)^{UV}\nabla_V\omega(U) - \omega([U, V]) \tag{82}$$

A simple computation shows that the exterior derivation of a simple tensor $\tilde{\omega} \otimes X$, is as follows.

$$d^\nabla(\tilde{\omega} \otimes X) = d\tilde{\omega} \otimes X + (-1)^{\tilde{\omega}}\tilde{\omega} \wedge \nabla X \tag{83}$$

In above, $\nabla X$ is the $E$-valued differential form of order one that is defined by $(\nabla X)(U) = \nabla_U X$. For a section $X$ of $E$, as an $E$-valued zero differential form we have $d^\nabla(X) = \nabla X$ and $d^\nabla d^\nabla(X)(U, V) = R(U, V)(X)$. The form $\nabla X$ is sum of two differential forms of pair orders (1,0) and (0,1). We can decompose $\nabla X$ as follows.

$$\nabla^{(1,0)}X = \sum_i E^i \wedge \nabla_{E_i}X \ , \ \nabla^{(0,1)}X = \sum_j \alpha^j \wedge \nabla_{\alpha_j}X \tag{84}$$

Clearly, $\nabla X = \nabla^{(1,0)}X + \nabla^{(0,1)}X$ and the pair order of $\nabla^{(1,0)}X$ is $(1,0)$ and the pair order of $\nabla^{(1,0)}X$ is $(0,1)$. In general, for a vector valued differential form $\omega$ of the pair order $(k, l)$, $d^\nabla \omega$ is sum of two differential forms of the pair orders $(k+1, l)$ and $(k, l+1)$. So, we can decompose the operator $d^\nabla$ to find two operators $d_1^\nabla$ and $d_2^\nabla$ that are defined as follows.

$$d_1^\nabla(\tilde{\omega} \otimes X) = d\tilde{\omega} \otimes X + (-1)^{\tilde{\omega}} \tilde{\omega} \wedge \nabla^{(1,0)}X \tag{85}$$

$$d_2^\nabla(\tilde{\omega} \otimes X) = \tilde{\omega} \wedge \nabla^{(0,1)}X \tag{86}$$

Clearly, $d^\nabla = d_1^\nabla + (-1)^{\tilde{\omega}} d_2^\nabla$, and if the pair order of $\omega$ is $(k, l)$, then the pair order of $d_1^\nabla \omega$ is $(k+1, l)$ and the pair order of $d_2^\nabla \omega$ is $(k, l+1)$.

If we choose a torsion free super connection on $M$, then the covariant derivation of an $E$-valued super differential form can be fined as follows.

$$\nabla_U(\tilde{\omega} \otimes X) = \nabla_U \tilde{\omega} \otimes X + \tilde{\omega} \otimes \nabla_U X \tag{87}$$

**Theorem (7.1):** For an $E$-valued super differential form $\omega$, the following relation holds.

$$(d_1^\nabla - d_2^\nabla)\omega = \sum E^i \wedge \nabla_{E_i} \omega - \sum \alpha^j \wedge \nabla_{\alpha_j} \omega \tag{88}$$

**Proof**: We prove this relation by the following computations.

$$\begin{aligned}
&\sum_i E^i \wedge \nabla_{E_i}(\tilde{\omega} \otimes X) - \sum_j \alpha^j \wedge \nabla_{\alpha_j}(\tilde{\omega} \otimes X) \\
&= \sum_i E^i \wedge \left((\nabla_{E_i}\tilde{\omega}) \otimes X + \tilde{\omega} \otimes \nabla_{E_i} X\right) - \sum_j \alpha^j \wedge \left((\nabla_{\alpha_j}\tilde{\omega}) \otimes X + \tilde{\omega} \otimes \nabla_{\alpha_j} X\right) \\
&= \left(\sum_i E^i \wedge (\nabla_{E_i}\tilde{\omega}) - \sum_j \alpha^j \wedge (\nabla_{\alpha_j}\tilde{\omega})\right) \otimes X \\
&\quad + (-1)^{\tilde{\omega}}\left(\tilde{\omega} \wedge \sum_i E^i \wedge \nabla_{E_i} X\right) - \sum_j \tilde{\omega} \wedge \alpha^j \otimes \nabla_{\alpha_j} X \\
&= d\tilde{\omega} \otimes X + (-1)^{\tilde{\omega}} \tilde{\omega} \wedge \nabla^{(1,0)}X - \sum_j \tilde{\omega} \wedge \nabla^{(0,1)}X \\
&= \left(d_1^\nabla - d_2^\nabla\right)(\tilde{\omega} \otimes X) \qquad \bullet
\end{aligned} \tag{89}$$

For an exterior product, associated to a product that is parallel with respect to super connections, by a simple computation, we obtain the following relation.

$$\begin{aligned}
d^\nabla(\omega \wedge \varphi) &= (d^\nabla \omega) \wedge \varphi + (-1)^\omega \omega \wedge d^\nabla \varphi \\
d_1^\nabla(\omega \wedge \varphi) &= (d_1^\nabla \omega) \wedge \varphi + (-1)^\omega \omega \wedge d_1^\nabla \varphi
\end{aligned} \tag{90}$$

Note that in above, different super bundles and super connections may occur.

## 8. Super metric on super bundles

In this section, we extend the notion of (semi)Riemannian metrics on vector bundles to super metrics on super vector bundles.

**Definition(8.1):** If $E$ is a super vector bundle on $M$, an even covariant super tensor of order two on $E$, $<.,.>: \Gamma E \times \Gamma E \to C^\infty(M)$ is called super metric, if for each $p \in M$ it is a super inner product on $E_p$.

So, even and odd sections are orthogonal to each other, and for homogeneous section $X$ and $Y$, we have:

$$< X, Y > = (-1)^{XY} < Y, X > \qquad (91)$$

A super connection $\nabla$ on $E$ is called compatible with a super metric $<.,.>$, whenever for each $U \in \widehat{\mathfrak{X}}M$ and $X, Y \in \Gamma E$, the following relation holds.

$$U < X, Y > = < \nabla_U X, Y > + (-1)^{UX} < X, \nabla_U Y > \qquad (92)$$

In the case $E = \widehat{T}M$, a super metric on $\widehat{T}M$ is called a super metric on $M$. The restriction of a super metric on $M$ to $TM$, is a semi-Riemannian metric on $M$. So, a super metric on $M$ is a pair $(g^e, g^o)$ in which $g^e$ is a semi-Riemannian metric on $M$ and $g^o$ is a nondegenerate symplectic form on the trivial bundle $M \times W$.

**Theorem(8.2):** If $<.,.>$ is a super metric on $M$, then there exists a unique torsion free super connection on $M$ that is compatible with this metric.

**Proof:** An extension of Koszul formula can be employed to define a torsion free super connection that is compatible with the metric. The following formula define a super connection that is torsion free and metric compatible.

$$\begin{aligned} 2 < \nabla_U V, W > = {} & U < V, W > + (-1)^{UV} V < U, W > - (-1)^{W(U+V)} W < U, V > \\ & + < [U, V], W > - < U, [V, W] > - (-1)^{VW} < [U, W], V > \end{aligned} \qquad (93)$$

It can easily be shown that the right hand expression is tensorial with respect to $W$, so by nondegenerancy of the super metric, the vector field $\nabla_U V$ is well defined. By the same method, we find that $\nabla$ is a super connection on $M$. Since, the inner product and Lie bracket preserve parities, $\nabla$ also preserve parities (i.e. $|\nabla_U V| = |U| + |V|$ ). To prove $\nabla$ is torsion free, we compute $\nabla_U V - (-1)^{UV} \nabla_V U$.

$$
\begin{aligned}
2 < \nabla_U V - (-1)^{UV} \nabla_V U, W > &= U < V, W > +(-1)^{UV} V < U, W > -(-1)^{W(U+V)} W < U, V > \\
&+ < [U,V], W > - < U, [V,W] > -(-1)^{VW} < [U,W], V > \\
&-(-1)^{UV} V < U, W > - U < V, W > +(-1)^{W(U+V)+UV} W < V, U > \\
&-(-1)^{UV} < [V,U], W > +(-1)^{UV} < V, [U,W] > +(-1)^{UW+UV} < [V,W], U > \\
&= 2 < [U,V], W >
\end{aligned} \quad (94)
$$

To prove $\nabla$ is metric compatible, we compute $< \nabla_U V, W > +(-1)^{UV} < V, \nabla_U W >$.

$$
\begin{aligned}
2(< \nabla_U V, W > +(-1)^{UV} < V, \nabla_U W >) &= 2(< \nabla_U V, W > +(-1)^{WV} < \nabla_U W, V >) \\
&= U < V, W > +(-1)^{UV} V < U, W > -(-1)^{W(U+V)} W < U, V > \\
&+ < [U,V], W > - < U, [V,W] > -(-1)^{VW} < [U,W], V > \\
&+(-1)^{WV} U < W, V > +(-1)^{UW+UV} W < U, V > -(-1)^{UV} V < U, W > \\
&+(-1)^{WV} < [U,W], V > -(-1)^{WV} < U, [W,V] > - < [U,V], W > \\
&= 2U < V, W >
\end{aligned} \quad (95)
$$

To prove uniqueness, assume $\nabla'$ is another torsion free super connection that is compatible with the metric. So,

$$
\begin{aligned}
2 < \nabla_U V, W > &= U < V, W > +(-1)^{UV} V < U, W > -(-1)^{W(U+V)} W < U, V > \\
&+ < [U,V], W > - < U, [V,W] > -(-1)^{VW} < [U,W], V > \\
&= < \nabla'_U V, W > +(-1)^{UV} < V, \nabla'_U W > +(-1)^{UV} < \nabla'_V U, W > + < U, \nabla'_V W > \\
&-(-1)^{W(U+V)} < \nabla'_W U, V > -(-1)^{WV} < U, \nabla'_W V > + < \nabla'_U V, W > +(-1)^{UV} < \nabla'_V U, W > \\
&- < U, \nabla'_V W > +(-1)^{WV} < U, \nabla'_W V > -(-1)^{VW} < \nabla'_U W, V > +(-1)^{VW+UW} < \nabla'_W U, V > \\
&= 2 < \nabla'_U V, W >
\end{aligned} \quad (96)
$$

This equality, implies $\nabla = \nabla'$. •

This connection is super version of the Levi-civita connection and it can be called Levi-civita super connection of the super metric.

**Theorem(8.3):** If $(g^e, g^o)$ is a super metric on $M$, and $\nabla = (\dot{\nabla}, \omega_0, X_{\alpha\beta})$ is its Levi-civita super connection, then $\dot{\nabla}$ is the Levi-civita connection of $g^e$ and

$$< \omega_0(U,\alpha), \beta > = < U, X_{\beta\alpha} > = \frac{1}{2} U(g^o)(\alpha, \beta) \quad (97)$$

**Proof**: The Koszul formula shows that $\dot{\nabla}$ is the ordinary Levi-civita connection of $g^e$, because for even vector fields, the extended Koszul formula reduces to the ordinary Koszul formula for ordinary metric $g^e$. The first equality is obtained from the following computation.

$$0 = \alpha <U, \beta> = <\nabla_\alpha U, \beta> + <U, \nabla_\alpha \beta> = <\omega_0(U, \alpha), \beta> + <U, X_{\alpha\beta}>$$
$$\Rightarrow <\omega_0(U, \alpha), \beta> = <U, X_{\beta\alpha}> \qquad (98)$$

The second equality can be obtained from the Koszul formula as follows.

$$\begin{aligned} 2 <U, X_{\beta\alpha}> &= 2 <U, \nabla_\beta \alpha> = 2 <\nabla_\beta \alpha, U> \\ &= \beta <\alpha, U> + (-1)^{\alpha\beta} \alpha <\beta, U> - (-1)^{U(\alpha+\beta)} U <\beta, \alpha> \\ &\quad + <[\beta, \alpha], U> - <\beta, [\alpha, U]> - (-1)^{\alpha U} <[\beta, U], \alpha> \\ &= -U\big(g^o(\beta, \alpha)\big) + g^o\big(\beta, U(\alpha)\big) + g^o\big(U(\beta), \alpha\big) \\ &= U(g^o)(\alpha, \beta) \qquad \qquad \qquad \qquad \qquad \qquad \bullet \end{aligned} \qquad (99)$$

**Theorem(8.4)**: The curvature tensor of the super Levi-civita connection has the following properties.

$$\begin{aligned} <R(U,V)(Z), W> &= -(-1)^{ZW} <R(U,V)(W), Z> \\ <R(U,V)(Z), W> &= (-1)^{(U+V)(Z+W)} <R(Z,W)(U), V> \end{aligned} \qquad (100)$$

**Proof**: The proof is the same as for ordinary curvature of a semi-Riemannian metric.●

By super contraction of the curvature tensor, Ricci curvature and scalar curvature are defined. First, consider a homogeneous local base $\{E_1, \cdots, E_n, \alpha_1, \cdots, \alpha_m\}$ for $\hat{T}M$ and its reciprocal base $\{E^1, \cdots, E^n, \alpha^1, \cdots, \alpha^m\}$. To define super Ricci curvature, first define $R(U, V, Z, W) = <R(U,V)(Z), W>$. The super contraction of the tensor $R(U, V, Z, W)$ with respect to $V$ and $Z$ is called the super Ricci curvature. In other words:

$$\boldsymbol{Ric}(U, V) = \sum_i <R(U, E_i)(E^i), V> + \sum_j <R(U, \alpha_j)(\alpha^j), V> \qquad (101)$$

Properties of the curvature tensor, show that the Ricci tensor is super symmetric.

$$\begin{aligned}
Ric(U,V) &= \sum_i <R(U,E_i)(E^i),V> + \sum_j <R(U,\alpha_j)(\alpha^j),V> \\
&= \sum_i (-1)^{UV} <R(E^i,V)(U),E_i> + \sum_j (-1)^{(1+U)(1+V)} <R(\alpha^j,V)(U),\alpha_j> \\
&= \sum_i (-1)^{UV} <R(V,E^i)(E_i),U> - \sum_j (-1)^{UV} <R(V,\alpha^j)(\alpha_j),U> \\
&= (-1)^{UV} Ric(V,U)
\end{aligned} \qquad (102)$$

In this computation, note that the reciprocal base of $\{E^1,\cdots,E^n,\alpha^1,\cdots,\alpha^m\}$ is $\{E_1,\cdots,E_n,-\alpha_1,\cdots,-\alpha_m\}$.

Scalar curvature is defined by super contraction of the Ricci curvature. So,

$$R = \sum_i Ric(E_i,E^i) + \sum_j Ric(\alpha_j,\alpha^j) \qquad (103)$$

**Divergence of super differential forms:** Here we can extend the notion of divergence for super differential forms. To this end, we mimic the relation between ordinary divergence and ordinary exterior derivation and the Hodge operator. Since, for each $p \in M$ we have a super inner product on $\hat{T}_p M = T_p M \oplus W$, we can use the Hodge operator for super differential forms.

**Definition(8.5):** For a numerical valued super differential form $\tilde{\eta}$ of the pair order $(k,l)$, divergence of $\tilde{\eta}$ is defined as follows.

$$\delta\tilde{\eta} = (-1)^{n(k+1)+s+1} H \circ d \circ H(\tilde{\eta}) \qquad (104)$$

If $k = 0$, then $\delta\tilde{\eta} = 0$, and for $1 \le k$ the pair order of $\delta\tilde{\eta}$ is $(k-1,l)$. For a simple tensor $\omega \otimes \xi$, its divergence can be computed as follows. Here, $\{E_1,\cdots,E_n\}$ is a local base of $TM$ and $\{E^1,\cdots,E^n\}$ is its dual base. For simplicity set $\varepsilon = (-1)^{n(k+1)+s+1}$.

$$\begin{aligned}
\delta(\omega \otimes \xi) &= \varepsilon H\left(d\left(H(\omega \otimes \xi)\right)\right) = \varepsilon H\left(d(H(\omega) \otimes \xi)\right) \\
&= \varepsilon H(d(H(\omega)) \otimes \xi + (-1)^{n-k} H(\omega) \wedge d\xi) \\
&= \varepsilon H(d(H(\omega)) \otimes \xi) + \varepsilon(-1)^{n-k} H(H(\omega) \wedge \sum_i E^i \wedge E_i(\xi)) \\
&= \delta\omega \otimes \xi + \varepsilon(-1)^{n-k} \sum_i H(H(\omega) \wedge E^i) \otimes E_i(\xi) \\
&= \delta\omega \otimes \xi - \sum_i i_{E^i}(\omega) \otimes E_i(\xi) = \delta\omega \otimes \xi - i_{d\xi}\omega
\end{aligned} \qquad (105)$$

In above, $i_{d\xi}\omega$ is an abbreviation of $\sum_i i_{E^i}(\omega) \otimes E_i(\xi)$ and it is an extension of insertion operator for super differential forms.

If $\xi$ is constant, then $d\xi = 0$ and $\delta(\omega\otimes\xi) = \delta\omega\otimes\xi$.

We can define point wise inner product of super differential forms. For any super differential forms $\tilde{\omega}$ and $\tilde{\varphi}$, $<\tilde{\omega},\tilde{\varphi}>$ is a smooth function on $M$ that is defined by $<\tilde{\omega},\tilde{\varphi}>(p) = <\tilde{\omega}_p,\tilde{\varphi}_p>$.

**Theorem(8.6):** If $\tilde{\omega}$ and $\tilde{\varphi}$ are super differential forms of the pair orders $(k,l)$ and $(k+1,l)$ respectively, and one of them is compactly supported, then

$$\int_M <d\tilde{\omega},\tilde{\varphi}>\Omega = \int_M <\tilde{\omega},\delta\tilde{\varphi}>\Omega \qquad (106)$$

**Proof:** A simple computation shows the equality.

$$\int <d\tilde{\omega},\tilde{\varphi}>\Omega = \int d\tilde{\omega}\wedge H(\tilde{\varphi}) = \int d(\tilde{\omega}\wedge H(\tilde{\eta})) - (-1)^k \tilde{\omega}\wedge dH(\tilde{\eta})$$
$$= \int -(-1)^{k+k(n-k)+s} \tilde{\omega}\wedge HHdH(\tilde{\eta}) \qquad (107)$$
$$= \int \tilde{\omega}\wedge H(\delta\tilde{\eta}) = \int <\tilde{\omega},\delta\tilde{\eta}>\Omega \quad \bullet$$

**Theorem(8.7):** For a homogeneous local base $\{E_1,\cdots,E_n,\alpha_1,\cdots,\alpha_m\}$, if $\tilde{\eta}$ is a super differential form, then $\delta\tilde{\eta}$ can be computed as follows.

$$\delta\tilde{\eta} = -\sum_i i_{E^i}\left(\nabla_{E_i}\tilde{\eta}\right) + \sum_j \left(\nabla_{\alpha_j} i_{\alpha^j}\tilde{\eta}\right) \qquad (108)$$

**Proof:** We prove the theorem for a simple tensor $\omega\otimes\xi$ of the pair order $(k,l)$. We can assume that $\xi$ is constant. So, $\delta(\omega\otimes\xi) = \delta\omega\otimes\xi$.

First, note that for $\xi \in C^\infty(M,\vee^l W^*)$, we have $\tilde{\omega} = \sum_j \left(\nabla_{\alpha^j} i_{\alpha_j}\xi\right) = 0$. Because, the pair order of $\tilde{\omega}$ is $(1,l-2)$. For $l \le 1$, clearly $\tilde{\omega} = 0$. For $2 \le l$ the following computations prove the equality.

$$\tilde{\omega}(V,\beta_1,\cdots,\beta_{l-2}) = \left(\sum_j \left(\nabla_{\alpha^j} i_{\alpha_j}\xi\right)\right)(V,\beta_1,\cdots,\beta_{l-2})$$
$$= \sum_j i_{\alpha_j}\xi(\nabla_{\alpha^j}V,\beta_1,\cdots,\beta_{l-2}) = \sum_j \xi(\alpha_j,\nabla_{\alpha^j}V,\beta_1,\cdots,\beta_{l-2})$$
$$= \sum_j \xi(<\alpha^t,\nabla_{\alpha^j}V>\alpha_t,\alpha_j,\beta_1,\cdots,\beta_{l-2}) \qquad (109)$$
$$= \sum_j <\alpha^t,\nabla_{\alpha^j}V>\xi(\alpha_t,\alpha_j,\beta_1,\cdots,\beta_{l-2})$$
$$= \sum_j <\nabla_{\alpha^j}\alpha^t,V>\xi(\alpha_t,\alpha_j,\beta_1,\cdots,\beta_{l-2}) = 0$$

Now, we compute the right hand side of the equality.

$$-\sum_i i_{E^i} \nabla_{E_i}(\omega \otimes \xi) + \sum_j \nabla_{\alpha_j} i_{\alpha^j}(\omega \otimes \xi) = -\sum_i i_{E^i} \nabla_{E_i}(\omega \otimes \xi) + (-1)^k \sum_j \nabla_{\alpha_j}(\omega \otimes i_{\alpha^j}\xi)$$

$$= -\sum_i i_{E^i}\left((\nabla_{E_i}\omega)\otimes\xi + \omega\otimes\nabla_{E_i}\xi\right) + (-1)^k \sum_j \left((\nabla_{\alpha_j}\omega)\wedge i_{\alpha^j}\xi + (-1)^k \omega \wedge \nabla_{\alpha_j} i_{\alpha^j}\xi\right)$$

$$= -\sum_i i_{E^i}(\nabla_{E_i}\omega)\otimes\xi - \sum_i i_{E^i}\omega \otimes \nabla_{E_i}\xi + (-1)^k \sum_j (\nabla_{\alpha_j}\omega)\wedge i_{\alpha^j}\xi$$

$$= \delta\omega\otimes\xi - \left(\sum_i i_{E^i}\omega \otimes \nabla_{E_i}\xi - (-1)^k \sum_j (\nabla_{\alpha_j}\omega)\wedge i_{\alpha^j}\xi\right) \qquad (110)$$

We must prove the second term is zero. The second term is a differential form of the pair order $(k-1, l)$.

$$\left(\sum_i i_{E^i}\omega \otimes \nabla_{E_i}\xi - (-1)^k \sum_j (\nabla_{\alpha_j}\omega)\wedge i_{\alpha^j}\xi\right)(V_1, \cdots, V_{k-1}, \beta_1, \cdots, \beta_l)$$

$$= \sum_i \omega(E^i, V_1, \cdots, V_{k-1})(\nabla_{E_i}\xi)(\beta_1, \cdots, \beta_l)$$

$$-(-1)^k \sum_{j,t} \omega\left(V_1, \cdots, V_{k-1}, \nabla_{\alpha_j}\beta_t\right)\xi(\alpha^j, \beta_1, \cdots, \hat{\beta}_t, \cdots, \beta_l)$$

$$= \sum_i \omega(E^i, V_1, \cdots, V_{k-1}) E_i\left(\xi(\beta_1, \cdots, \beta_l)\right) - \sum_{i,t} \omega(E^i, V_1, \cdots, V_{k-1})\xi(\beta_1, \cdots, \nabla_{E_i}\beta_t, \cdots, \beta_l)$$

$$-(-1)^k \sum_{i,j,t} \omega\left(V_1, \cdots, V_{k-1}, <E_i, \nabla_{\alpha_j}\beta_t> E^i\right)\xi(\alpha^j, \beta_1, \cdots, \hat{\beta}_t, \cdots \beta_l) \qquad (111)$$

$$= \sum_{i,t} \omega(E^i, V_1, \cdots, V_{k-1})\xi\left(\beta_1, \cdots, E_i(\beta_t) - \nabla_{E_i}\beta_t, \cdots, \beta_l\right)$$

$$-(-1)^k \sum_{i,j,t} \omega(V_1, \cdots, V_{k-1}, <\omega_0(E_i, \beta_t), \alpha_j> E^i)\xi(\alpha^j, \beta_1, \cdots, \hat{\beta}_t, \cdots \beta_l)$$

$$= -\sum_{i,t} \omega(E^i, V_1, \cdots, V_{k-1})\xi(\beta_1, \cdots, \omega_0(E_i, \beta_t), \cdots \beta_l)$$

$$-(-1)^k \sum_{i,t} \omega(V_1, \cdots, V_{k-1}, E^i)\xi(\omega_0(E_i, \beta_t), \beta_1, \cdots, \hat{\beta}_t, \cdots \beta_l) = 0 \qquad \bullet$$

For a super vector bundle $E$, that equipped with a super connection, the notion of divergence is defined for $E$-valued super differential forms, in the same way. But, in this case the exterior derivation $d^\nabla$ is decomposed to two operator $d_1^\nabla$ and $d_2^\nabla$ and it is better to define divergence associate to each of these operators separately. The divergence operator associated to $d_1^\nabla$ is denoted by $\delta_1^\nabla$ and is defined as follows.

$$\delta_1^\nabla \varphi = (-1)^{n(k+1)+s+1} H \circ d_1^\nabla \circ H(\varphi) \qquad (112)$$

If the pair order of $\tilde{\varphi}$ is $(k, l)$ then the pair order of $\delta_1^\nabla \tilde{\varphi}$ is $(k-1, l)$. For a simple tensor $\tilde{\omega}\otimes X$ a similar computation as equalities in (105) shows that

$$\delta_1^\nabla(\tilde{\omega}\otimes X) = \delta\tilde{\omega}\otimes X - \sum_i (i_{E^i}\tilde{\omega})\otimes \nabla_{E_i} X \qquad (113)$$

If $E$ has a super metric that $\nabla$ is compatible with respect to it, a function valued product of $E$-valued super differential forms can be defined as follows.

$$< \tilde{\omega} \otimes X, \tilde{\varphi} \otimes Y > = < \tilde{\omega}, \tilde{\varphi} > < X, Y > \qquad (114)$$

**Theorem(8.8):** If $\omega$ and $\varphi$ are $E$-valued super differential forms of the pair orders $(k, l)$ and $(k + 1, l)$, respectively, and one of them is compactly supported, then

$$\int_M < d_1^\nabla \omega, \eta > \Omega = \int_M < \omega, \delta_1^\nabla \eta > \Omega \qquad (119)$$

**Proof:** Computation for this theorem is the same as theorem (8.6). ●

We define the divergence operator associated to $d_2^\nabla$ such that equality of the theorem (8.8) satisfied. To this end, define $\delta_2^\nabla$ as follows.

$$\delta_2^\nabla(\tilde{\omega} \otimes X) = i_{\alpha^j} \tilde{\omega} \otimes \nabla_{\alpha_j} X \qquad (120)$$

**Theorem(8.9):** The following relation holds between $d_2^\nabla$ and $\delta_2^\nabla$.

$$< d_2^\nabla(\tilde{\omega} \otimes X), \tilde{\varphi} \otimes Y > = (-1)^X < \tilde{\omega} \otimes X, \delta_2^\nabla(\tilde{\varphi} \otimes Y) > \qquad (121)$$

**Proof:** Assume $\tilde{\omega} = \omega \otimes \xi$ and $\tilde{\varphi} = \varphi \otimes \eta$

$$\begin{aligned}
< d_2^\nabla(\tilde{\omega} \otimes X), \tilde{\varphi} \otimes Y > &= < \sum_j \omega \otimes (\xi \vee \alpha^j) \otimes \nabla_{\alpha_j} X, \varphi \otimes \eta \otimes Y > \\
&= \sum_j < \omega, \varphi > < \xi \vee \alpha^j, \eta > < \nabla_{\alpha_j} X, Y > \\
&= \sum_j < \omega, \varphi > \left(-< \xi, i_{\alpha^j} \eta >\right)\left(-(-1)^X < X, \nabla_{\alpha_j} Y >\right) \qquad (122)\\
&= \sum_j (-1)^X < \omega \otimes \xi \otimes X, \varphi \otimes i_{\alpha^j} \eta \otimes \nabla_{\alpha_j} Y > \\
&= (-1)^X < \tilde{\omega} \otimes X, \delta_2^\nabla(\tilde{\varphi} \otimes Y) >
\end{aligned}$$

**Theorem(8.10):** For a torsion free super connection on $M$, and a homogenous local base $\{E_1, \cdots, E_n, \alpha_1, \cdots, \alpha_m\}$, the operator $\delta^\nabla = \delta_1^\nabla + \delta_2^\nabla$ can be computed as follows.

$$\delta^\nabla \omega = -\sum_i i_{E^i}\left(\nabla_{E_i} \omega\right) + \sum_j \left(\nabla_{\alpha_j} i_{\alpha^j} \omega\right) \qquad (123)$$

**Proof:** We compute right hand side of the equality given above for $\omega = \tilde{\omega} \otimes X$.

$$\begin{aligned}
&-\sum_i i_{E^i}\nabla_{E_i}(\tilde{\omega}\otimes X) + \sum_j \nabla_{\alpha_j}i_{\alpha^j}(\tilde{\omega}\otimes X)\\
&= -\sum_i i_{E^i}\left(\left(\nabla_{E_i}\tilde{\omega}\right)\otimes X + \tilde{\omega}\otimes\nabla_{E_i}X\right) + \sum_j \left(\left(\nabla_{\alpha_j}i_{\alpha^j}\tilde{\omega}\right)\otimes X + i_{\alpha^j}\tilde{\omega}\wedge\nabla_{\alpha_j}X\right)\\
&= -\sum_i i_{E^i}\left(\nabla_{E_i}\tilde{\omega}\right)\otimes X + \sum_j \left(\nabla_{\alpha_j}i_{\alpha^j}\tilde{\omega}\right)\otimes X - \sum_i i_{E^i}\tilde{\omega}\otimes\nabla_{E_i}X + \sum_j i_{\alpha^j}\tilde{\omega}\wedge\nabla_{\alpha_j}X\\
&= \delta\tilde{\omega}\otimes X - \sum_i i_{E^i}\tilde{\omega}\otimes\nabla_{E_i}X + \sum_j i_{\alpha^j}\tilde{\omega}\wedge\nabla_{\alpha_j}X = \delta_1^{\nabla} + \delta_2^{\nabla} = \delta^{\nabla} \qquad \bullet
\end{aligned} \qquad (124)$$